\documentclass{article}
\usepackage{sw20bams}
\usepackage{latexsym}
\usepackage{amsmath}
\usepackage{amsfonts}
\usepackage{amssymb}
\usepackage{graphicx}
\newtheorem{theorem}{Theorem}
\newtheorem{acknowledgement}[theorem]{Acknowledgement}

\newtheorem{definition}[theorem]{Definition}

\newtheorem{lemma}[theorem]{Lemma}
\newtheorem{notation}[theorem]{Notation}

\newtheorem{proposition}[theorem]{Proposition}
\newtheorem{remark}[theorem]{Remark}

\begin{document}

\author{Jay Jorgenson\\City College of NY\\Department of Mathematics,\\New York, NY 10031.
\and Andrey Todorov\\University of California\\Department of Mathematics\\Santa Cruz, CA 95064\\Bulgarian Academy of Sciences\\Institute of Mathematics\\Ul. G. Bonchev 8\\Sofia, Bulgaria.}
\title{Ample Divisors, Automorphic Forms and Shafarevich's Conjecture}
\date{June 4/1999}
\maketitle
\begin{abstract}
In this article we give a general approach to the following analogue of
Shafarevich's conjecture for some polarized algebraic varieties; suppose that
we fix a type of an algebraic variety and look at families of such type of
varieties over a fixed Riemann surface with fixed points over which we have
singular varieties, then one can ask if the set of such families, up to
isomorphism, is finite.

In this paper we give a general approach to such types of problems. The main
observation is the following; suppose that the moduli space of a fixed type of
algebraic polarized variety exists and suppose that in some projective smooth
compactification of the coarse moduli the discriminant divisor supports an
ample one, then it is not difficult to see that this fact implies the analogue
of Shafarevich's conjecture.

In this article we apply this method to certain polarized algebraic K3
surfaces and also to Enriques surfaces.
\end{abstract}
\tableofcontents

\section{Introduction}

\subsection{General Historical Review of the Problem}

In his article published in the Proceedings of the International Congress of
Mathematicians, Stockholm meeting held in 1962, Shafarevich wrote:

\textit{''One of the main theorems on algebraic numbers connected with the
concept of discriminant is Hermit's theorem, which states that the number of
extensions }$k^{\prime}/k$ \textit{of a given degree and given discriminant is
finite. This theorem may be formulated as follows: the number of extensions
}$k^{\prime}/k$\textit{\ of a given degree whose critical prime divisors
belong to a given finite set S is finite.''}

Inspired by this result of Hermit, Shafarevich conjectured: \textit{''There
exists only a finite number of fields of algebraic functions }$K/k$%
\textit{\ of a given genus }$g\geq1$\textit{, the critical prime divisors of
which belong to a given finite set }$S."$ See page 290 of \cite{sh}.

In unpublished work, Shafarevich proved his conjecture in the setting of
hyperelliptic curves. On page 755 of \cite{sh}, in his remarks on his papers,
he wrote: \textit{''Here two statements made in lecture are mixed into one:
formulation of a result and conjecture. The result was restricted to the case
of a hyperelliptic curves while conjecture concerned general curves...This
conjecture became much more attractive after A. N. Parshin proved that implies
Mordell conjecture...In 1983 it was proved by G. Faltings (Invent. Math.
\textbf{73, }}439-466(1983))\textit{, with a proof of the Mordell conjecture
as a consequence.''}

In this paper we developed a general approach \ to higher dimensional
analogues of Shafarevich's type conjectures over function field and apply this
method to families of algebraic polarized K3 surfaces of special degrees or to
families of polarized Enriques surfaces over an algebraic curve.

The formulation of the problem is the following. Let C be a fixed curve of
genus g and let E be a set of different points on C. Let C be a fixed,
non-singular algebraic curve, and let E be a fixed effective divisor on C such
that all the points in E have multiplicity 1. Define Sh(C,E,Z) to be the set
of all isomorphism classes of projective algebraic varieties
$\mathcal{Z\rightarrow}C$ with a fibre of ''type' $Z$ such that the singular
fibres are over the set E. The general Shafarevich type problem
is:\textbf{\ ''}\textit{For which ''type'' of varieties }$Z$ \textit{and data
(C,E)\ is such that }Sh (C,E,Z) \textit{is finite} ?

Previous work on Shafarevich-type problems include the following results. In
the case when Z is a curve of genus $g>1$ and E is empty, Parshin proved that
Sh (C,E,Z) is finite, and Arakelov proved finiteness in the case E is not
empty. See \cite{Ar} and \cite{Par}. Faltings constructed examples showing
that Sh (C,E,Z) is infinite for abelian varieties of dimension $\geq8.$ See
\cite{Falt}. Saito an Zucker extended the construction of Faltings to the
setting when Z is an algebraic polarized K3 surface. They were able to
classify all cases when the set Sh (C,E,Z) is infinite. They are not
considering polarized families. See \cite{SZ}. Faltings proved the
Shafarevich's conjecture over the number fields and thus he proved Mordell
conjecture. Yves Andre proved the analogue of Shafarevich's conjectures over
the number fields for K3 surfaces. See \cite{an}. Using techniques from
harmonic maps Jost and Yau analyzed Sh (C,E,Z) for a large class of varieties.
See \cite{JY}. Ch. Peters studied finiteness theorems by considering
variations of Hodge structures and utilizing differential geometric aspects of
the period map and associated metrics on the period domain. See \cite{Pe}.

Our approach to prove finiteness of Sh (C,E,Z) is to analyze the discriminant
locus $\mathcal{D}$ in some compactification of the coarse moduli space$,$
i.e. those points which correspond to singular varieties in the Baily Borel
compactification of the coarse moduli space of pseudo polarized algebraic K3
surfaces or Enriques surfaces. The main observation is that if $\mathcal{D}$
supports an ample divisor in some compactification, then Sh (C,E,Z) is finite.

The second author was informed that A. Parshin and E. Bedulev have proved
finiteness for family of algebraic surfaces over a fixed algebraic curve
assuming that all the fibres are non-singular.

Recently Migliorni, Kov\'{a}cs and Zhang that any family of minimal algebraic
surfaces of general type over a curve of genus g and m singular points such
that 2g-2+m$\leq0$ is isotrivial. This result was recently reproved by E.
Bedulev and E. Viehweg. See \cite{Ko}, \cite{Ko1}, \cite{Ko2}, \cite{Ko3},
\cite{Mi}, \cite{OW}, \cite{WZ}, \cite{Z} and \cite{BW}.

\subsection{Organization of the Article}

In \textbf{Section 2}$,$ we define the terms needed to properly formulate the
problem, and we prove general rigidity result from which the corresponding
finiteness will follow. The main observation is that if the discriminant locus
$\mathcal{D}$ in some compactification of the course moduli space of certain
type of polarized algebraic variety $\mathcal{M}$ is an ample divisor, then
the moduli space $\frak{M}(C;p_{1},..,p_{k})$ of maps of a pair of a fixed
algebraic curve with fixed points $p_{1},..,p_{k}$ on it $(C;p_{1,..,}p_{k})$
to $\mathcal{M}$ such that the points $p_{1},..,p_{k} $ are mapped to
$\mathcal{D},$ then $\frak{M}(C;p_{1},..,p_{k})$ is a discrete set. It is a
well known fact that in order to prove finiteness of $\frak{M}(C;p_{1}%
,..,p_{k}),$ one needs to prove that the volume of the image of $C$ in
$\mathcal{M}$ is bounded with respect to some metric on $\mathcal{M}.$

In \textbf{Section 3} we discuss various background material in the study of
K3 surfaces following \cite{Ast} and \cite{JT}.

In \textbf{Section 4} we discuss various background material in the study
Enriques surfaces following \cite{JT98} and \cite{JT99}.

In \textbf{Section 5} we study the question when the discriminant locus in
some compactification of the moduli space of polarized algebraic K3 surfaces
is supporting an ample divisor. According to the Torelli Theorem and the
epimorphism of the period map we know that the moduli space of algebraic
polarized K3 surfaces is a locally symmetric space $\Gamma_{K3,2d}%
\backslash\frak{H}_{2,19},$ where $\frak{H}_{2,19}=SO_{0}(2,19)/SO(2)\times
SO(19)$ and $\Gamma_{K3,2d}$ is an arithmetic group acting on $\frak{H}%
_{2,19}.$ In \cite{JT} we proved that if Baily-Borel compactification
$\overline{\Gamma_{K3,2d}\backslash\frak{H}_{2,19}}$ of $\Gamma_{K3,2d}%
\backslash\frak{H}_{2,19}$ contains only one cusp of dimension $0,$ then there
exists an automorphic form $\eta_{2d}$ such that the support of the zero set
of $\eta_{2d}$ is exactly the closure $\overline{\mathcal{D}}$ of the
discriminant locus $\mathcal{D}.\footnote{\textbf{1. }In \cite{JT} we refer to
a Theorem of F. Scattone which is not correct. \textbf{2.} It was Nikulin who
pointed out a mistake in an earlier version of \cite{JT}.}$ Horikawa and the
second author proved that for algebraic K3 surfaces with a polarization class
$e$ such that $\left\langle e,e\right\rangle =2$ the Baily-Borel
compactification $\overline{\Gamma_K3,2d\backslash\frak{H}_2,19}$ contains
only one cusp of dimension 0. See \cite{Hor}. Borcherds constructed in
\cite{BKPS} a holomorphic automorphic form in case of degree two polarization
based on idea of the second author.

In \textbf{Section 6} by using Gauss-Bonnet Theorem derive the following two Theorems:

\textbf{THEOREM.} Let C be an algebraic curve of genus g and $\pi:Y\rightarrow
C$ be a three dimensional projective non-singular variety such that for every
$t\in C\backslash E,$ $\pi^{-1}(t)=X_{t}$ is a non-singular K3 surface and for
each $t\in E,$ $\pi^{-1}(t)=X_{t}$ is a singular surface. Suppose that on $Y$
we have a polarization class $H$ such that $H\left|  _{X_{t}}\right.  =e,$
$<$%
e,e%
$>$%
=2d and the Baily-Borel compactification $\overline{\Gamma_{K3,2d}%
\backslash\frak{H}_{2,19}}$ contains only one cusp of dimension 0. Let
m$_{\infty}$ be the number of points on C for which the local monodromy
operator is of infinite order. Then the number of singular fibres of $\pi$ is
less or equal to $2g-2+m_{\infty}.$

\textbf{THEOREM.} Let C be an algebraic curve of genus g and $\pi:Y\rightarrow
C$ be a three dimensional projective non-singular variety such that for every
$t\in C\backslash E,$ $\pi^{-1}(t)=X_{t}$ is a non-singular Enriques surface
and for each $t\in E,$ $\pi^{-1}(t)=X_{t}$ is a singular surface. Let
m$_{\infty}$ be the number of points on C for which the local monodromy
operator is of infinite order. Then the number of singular fibres of $\pi$ is
less or equal to $2g-2+m_{\infty}.$

In \textbf{Section 7 }by using the results of \textbf{Section 2} and the
Theorem just formulated above we prove the finiteness of Sh (C,E,\{K3,2\}). By
using the results of Saito and Zucker we prove the analogue of Shafarevich's
conjecture for polarizations $e$ such that $\left\langle e,e\right\rangle =2d$
and the Baily-Borel compactification $\overline{\Gamma_{K3,2d}\backslash
\frak{H}_{2,19}}$ contains only one cusp of dimension 0.

In \textbf{Section 8} we prove that Sh (C,E,Z) is finite when Z is an S-K3
surface, which means that for each $t\in C\backslash E,$ $\pi^{-1}(t)=X_{t}$
is a K3 surface whose Picard Group$\subseteq S,$ where S is a special
primitive sublattice in $H^{2}(X,\mathbb{Z})$ of signature (1,k).

In \textbf{Section 9} we prove that Sh (C,E,Z) is finite when Z is an Enriques surface.

\begin{acknowledgement}
We thank B. Grauder, S. Lang and G. Zuckerman for their interest in this
article. The second author is grateful to M. Schneider and Y. T. Siu for the
invitation to work at MSRI. We want to thank S.-T. Yau for drawing our
attention to \cite{JY}. Both authors acknowledge conversations with R.
Borcherds, L. Katzarkov and T. Pantev on mathematical topics related to this
paper. See \cite{BKPS}. We want to thank L. Katzarkov for his help with
\textbf{Section 2}. The second author wants to thanks R. Borcherds for his
valuable comments about an earlier version of this paper.
\end{acknowledgement}

\section{Basic Rigidity Result}

\subsection{Introduction}

In this section we shall prove a general Theorem about rigidity from which we
shall derive finiteness results. The results in this section are established
in very general context.

\subsection{Basic Definitions and Facts about Moduli of Maps}

Let S and X be two projective varieties. Let $f:S\rightarrow X$ be a morphism
between them. Let $\Gamma_{f}\subset S\times X$ be the graph of the map
$f:S\rightarrow X$ . According to the results of Grothendieck the Hilbert
scheme of $\Gamma_{f}\subset S\times X$ will be a projective scheme. See
\cite{SGA}.

\begin{definition}
\label{Gro}We will denote the Hilbert scheme of $\Gamma_{f}\subset S\times X $
by $\frak{M}_{f}(S,X)$ and called the moduli space of the map $f.$
\end{definition}

Grothendieck proved that $\frak{M}_{f}(S,X)$ is a projective scheme. See
\cite{SGA}.

\begin{definition}
\label{d1}\textbf{a.} Let f:S$\rightarrow$X be a morphism of projective
varieties with dim S$\leq$dim X and suppose that the morphism f:S$\rightarrow
$f(S) is finite. We say that f admits a non trivial one parameter deformation
if there is a non-singular projective curve T and a family of algebraic maps
F:T$\times$S$\rightarrow$X such that for some t$_{0}\in$T we have F$_{t_{0}}%
=$f and the morphism F:T$\times$S$\rightarrow$F((T$\times$S) is finite too.
\end{definition}

\textbf{b. }\textit{We will say that the deformation of f is trivial if
F}$_{t}$=\textit{f for all t}$\in$T.

\textbf{c.}\textit{\ Two families of maps F}$_{1}:T\times$\textit{S}%
$_{1}\rightarrow$\textit{X and F}$_{2}:T\times$\textit{S}$_{2}\rightarrow
$\textit{X are said to be isomorphic if there is a common finite cover S of
S}$_{1}$\textit{\ and S}$_{2}$\textit{\ such that the lifts of F}$_{1}%
$\textit{\ to T}$\times$\textit{S and F}$_{2}$ \textit{to T}$\times$\textit{S
are isomorphic, meaning there exists a biholomorphic map id}$\times$\textit{g
from T}$\times$\textit{S to itself such that F}$_{1}$\textit{=F}$_{2}%
\circ(id\times\mathit{g)}$\textit{.}

\begin{remark}
From now on we will consider only one parameter deformations of maps.
\end{remark}

\subsection{The Rigidity Theorem}

\begin{lemma}
\label{l1}Let f:S$\rightarrow$X be a morphism between projective varieties S
and X such that f:S$\rightarrow$f(S) is a finite morphism, and let F be a
deformation of f. Let D be an ample Cartier divisor on X, and assume that the
image of S is not contained in D. Suppose that F$^{\ast}($D)$=$D$_{S}\times$T
where D$_{S}$ is a Cartier divisor on S, then F is the trivial deformation.
\end{lemma}

\textbf{PROOF:\ }Choose sufficiently large n such that nD is a very ample
divisor on X. Since Definition \ref{d1} implies that the map F is a finite map
we can conclude that F$^{\ast}$(nD)=nF$^{\ast}$(D) will be a very ample
divisor on S$\times$T. From the condition that F$^{\ast}($D)$=$D$_{S}\times$T
we deduce that the line bundle $\mathcal{O}_{\text{S}\times\text{T}}$
(F$^{\ast}($nD)) is isomorphic to $\mathcal{O}_{\text{T}}\otimes
\mathcal{O}_{\text{S}}(n$D$_{S}).$ So we have

\begin{center}
$H^{0}($S$\times$T,$\mathcal{O}_{\text{S}\times\text{T}}$ (F$^{\ast}%
($nD)))=$H^{0}($S$\times$T,$\mathcal{O}_{\text{T}}\otimes\mathcal{O}%
_{\text{S}}(n$D$_{S}))=H^{0}($T,$\mathcal{O}_{\text{T}})\otimes H^{0}%
($S,$\mathcal{O}_{\text{S}}(n$D$_{S})).$
\end{center}

From here we obtain that any section $\sigma\in H^{0}($S$\times$%
T,$\mathcal{O}_{\text{S}\times\text{T}}$ (F$^{\ast}($nD)))=$H^{0}%
($T,$\mathcal{O}_{\text{T}})\otimes H^{0}($S,$\mathcal{O}_{\text{S}}(n$%
D$_{S}))$ can be written in the form $\sigma=\sum_{i}c_{i}\otimes\pi_{2}%
^{\ast}(\sigma_{i})$ where $c_{i}$ are constants and $\sigma_{i}$ are sections
of $\mathcal{O}_{\text{S}}($D$_{S})$ on S. So we can conclude that the
projective morphism defined by the choice of the basis in $H^{0}($S$\times
$T,$\mathcal{O}_{\text{S}\times\text{T}}$ (F$^{\ast}($nD))) of the very ample
line bundle $\mathcal{O}_{\text{S}\times\text{T}}$ (F$^{\ast}($nD)) will map
the subvariety s$\times T$ \ to a point. This contradicts the definition of
the very ample line bundle, which implies that any basis in $H^{0}($S$\times
$T,$\mathcal{O}_{\text{S}\times\text{T}}$ (F$^{\ast}($nD))) will define an
embedding of S$\times$T into some projective variety. $\blacksquare.$

\begin{theorem}
\label{t1}Let S and X be projective varieties, and let D be an ample Cartier
divisor on X. Let $\mathcal{F}$ be the set of all finite morphisms
f:S$\rightarrow$X such that \textbf{a.} f(S) is not contained in D,
\textbf{b.} the map f:S$\rightarrow$f(S) is finite \textbf{c.} f$^{-1}%
$(D)=D$_{S}$, where D$_{S}$ is a fixed divisor on S for all f$\in\mathcal{F}$.
\ \textbf{d. }Suppose that vol(f(S))%
$<$%
$c,$ where the volume of f(S) is with respect to the restriction of Fubini
Study metric on f(S) obtained from the embedding of X in some projective space
given by the very ample divisor $n_{0}$D, then $\mathcal{F}$ \ is a finite set.
\end{theorem}

\textbf{PROOF:} Let us choose n sufficiently parge positive integer such that
nD is a very ample divisor on X. By assumption, the divisor $\mathcal{D}%
=n$f$^{\ast}($D)$\times$X+S$\times$nD is an ample divisor on S$\times$X since
f is a finite. Let us denote by $\left|  \mathcal{D}\right|  $ the linear
system defined by the ample divisor $\mathcal{D}$. Let us fix a finite map
f:S$\rightarrow$X such that f fulfils conditions \textbf{a, b }and \textbf{c. }

\begin{definition}
\label{f}Once we fix the map f as above we define the set $\mathcal{M}_{f,D}$
as the subset of deformation space $\frak{M}_{f}(S,X)$ of the fixed map f as
defined in \textbf{Definition \ref{Gro}}, which fulfils conditions \textbf{a,
b }and \textbf{c }of Theorem \ref{t1}.
\end{definition}

\begin{proposition}
\label{f1}$\mathcal{M}_{f,D\text{ }}$ is a finite set.
\end{proposition}

\textbf{PROOF}: Since f is a finite map it is easy to see that the linear
system $|\mathcal{D}|$ defines a finite projective map $\phi_{|\mathcal{D}|}%
:$S$\times$X$\rightarrow\mathbb{P}^{N}$ once we choose a basis in $H^{0}%
($S$\times$X,$\mathcal{O}_{\text{S}\times\text{X}}(\mathcal{D})).$ Let
$\Gamma_{\text{f}}$ $\subset$S$\times$X be the graph of the map
f:S$\rightarrow$X and let $\frak{M}_{f}(S,X)$ be the projective variety
defined in Definition \ref{Gro}. From the results in \cite{SGA} we deduce that
the condition that f(S) is not contained in the fixed Cartier divisor D define
a Zariski open set in\ $\frak{M}_{f}(S,X)$ which we will denote by
$\frak{M}_{f,D}$. On the other hand the condition that f$\in\frak{M}_{f,D}$
and f$^{\ast}($D) is a fixed divisor D$_{S}$ in S is a closed condition, i.e.
the set $\mathcal{M}_{f,D}:=\{$f $\in\mathcal{M}_{D}$%
$\vert$%
f(D)$\varsubsetneq$D and f$^{\ast}($D) is a fixed divisor D$_{S}$ in S\} is a
closed subscheme in $\frak{M}_{f,D}.$ Lemma \ref{l1} implies that the set
$\mathcal{M}_{f,D}$ is zero dimensional. Since $\mathcal{M}_{f,D}$ is a
quasi-projective, zero dimensional scheme, then as a set $\mathcal{M}%
_{f,D\text{ }}$ is a finite set. Proposition \ref{f1} is proved. $\blacksquare.$

It remains to show that the set $\mathcal{Z}_{D}$ is a finite set where
$\mathcal{Z}_{D}$ is the union of all $\mathcal{M}_{f,D}$ , over all maps f
which fulfil the conditions \textbf{a, b, c, }and \textbf{d }.

\begin{proposition}
\label{f2}$\mathcal{Z}_{D}$ is a finite set.
\end{proposition}

\textbf{PROOF}: Define a height on $\mathcal{Z}$ in the following manner. Let
$\phi_{|n_{0}D|}:X\subset\mathbb{P}^{m}$ be the embedding given by the linear
system $|n_{0}D|$ of the very ample divisor $n_{0}D$ on X. Let $\omega
_{X}:=\phi_{|n_{0}D|}^{\ast}(\omega_{\mathbb{P}^{m}}),$ where $\omega
_{\mathbb{P}^{m}}$ be the Fubini-Study form. Each point of $\mathcal{Z}$ is
represented by a map $f:S\rightarrow X$ which satisfies conditions \textbf{a,
b, }and \textbf{c }of Theorem \ref{t1}\textbf{.} Define the height function
$h$ on $\mathcal{Z}$ as follows:

\begin{center}
$h(f$)=$\int_{f(S)}\wedge^{k}(\omega_{\mathbb{P}^{m}})=vol(f(S))=(n_{0}%
)^{k}<D,..,D>\left|  _{f(S)}\right.  $
\end{center}

where k is the dimension of S. Since f is a finite morphism, then f(S) will
have also dimension k. Since $h$ is defined as an intersection number, $h$ is
integer valued, hence $h$ is a locally constant function, i.e. constant on
each connected component of $\mathcal{Z}.$ Condition \textbf{d} implies that
$h$ is a bounded function.

We will prove now that $\mathcal{Z}_{D}$ is a compact set. Since
$\mathcal{Z}_{D}$ is a discrete set, therefore $\mathcal{Z}_{D}$ will be
finite. We need to prove that from any sequence $\{f_{n}\}$ in $\mathcal{Z}$
there is a subsequence which converges weakly to an algebraic map, i.e. the
corresponding subsequence of images of S converges to an algebraic subvariety
of X. Bishop's theorem implies that. See page 292 of \cite{Bi} or page 321 of
\cite{Gro}. Proposition \ref{f2} is proved. $\blacksquare.$

From Proposition \ref{f2} Theorem \ref{t1} follows directly. $\blacksquare.$

\begin{remark}
Our method of proof of Shafarevich-type problems for varieties over function
fields utilizes the results of this section in the following manner. In order
to prove finiteness of Sh(C,E,Z), we let X be the course moduli space of
varieties of type Z. Let $\overline{X}$ be a compactification of X such that
D=$\overline{X}$%
$\backslash$%
X is a divisor of normal crossings. It is then necessary to show that D
supports an ample divisor on $\overline{X}.$ We then study maps of C into
$\overline{X}$ with the requirement that the subset of C which intersects D is
exactly E. Theorem \ref{t1} yields the desired finiteness result.
\end{remark}

\section{Basic Properties of K3}

We will review some basic properties of algebraic K3 surfaces. For a more
general and complete discussion, the reader is referred to \cite{Ast} and
\cite{BPV}.

\subsection{Definition of a K3 Surface}

A K3 surface is a compact, complex two dimensional manifold with the following
properties: \textbf{i. }There exists a non-zero holomorphic two form $\omega$
on X. \textbf{ii. }$H^{1}($X,$\mathcal{O}_{\text{X}})=0.$

\begin{remark}
For the purposes of this article, we will assume that all surfaces are
projective varieties.
\end{remark}

From the defining properties, one can prove that the canonical bundle on X is
trivial. In \cite{Ast} and \cite{BPV}, the following topological properties
are proved. The surface X is simply connected, and the homology group
$H_{2}(X,\mathbb{Z})$ is a torsion free abelian group of rank 22. The
intersection form $<$ , $>$ on $H_{2}(X,\mathbb{Z})$ has the properties:
\textbf{i. }$<$u,u$>=0$ mod$(2);$\textbf{ii.} $\det\left(  \left\langle
e_{i},e_{j}\right\rangle \right)  =-1$ and \textbf{iii. }the symmetric form $<
$ , $>$ has a signature $(3,19).$

Theorem \textbf{5} on page 54 of \cite{Se} implies that as an Euclidean
lattice $H_{2}(X,\mathbb{Z})$ is isomorphic to the K3 lattice $\Lambda_{K3},$ where

\begin{center}
$H_{2}(X,\mathbb{Z})\backsimeq\Lambda_{K3}=\mathbf{H}^{3}\oplus(-E_{8})^{2}$
\end{center}

with

\begin{center}
$\mathbf{H=}\left(
\begin{array}
[c]{ll}%
0 & 1\\
1 & 0
\end{array}
\right)  $
\end{center}

being the hyperbolic lattice. Let $\alpha=\{\alpha_{i}\}$ be a basis of
$H_{2}(X,\mathbb{Z})$ with intersection matrix $\Lambda_{K3}.$ The pair
$(X,\alpha)$ is called a marked K3 surface. Let $e\in H^{1,1}(X,\mathbb{R}%
)\cap H^{2}(X,\mathbb{Z})$ be the class of a hyperplane section, i.e. an ample
divisor. The triple $(X,\alpha,e)$ is called a marked, polarized K3 surface.
The degree of the polarization is an integer $2d$ such that $<e,e>=2d.$

\subsection{Moduli of Marked, Algebraic and Polarized K3 surfaces}

From \cite{PS} and \cite{Ku} we have that the moduli space of isomorphism
classes of marked, polarized, algebraic K3 surfaces of a fixed degree $2d$,
which we denote by $\mathcal{M}_{K3,mpa}^{2d},$ is equal to an open set in the
symmetric space $\frak{h}_{2,19}:=SO_{0}(2,19)/[SO(2)\times SO(19)].$ Let

\begin{center}
$\Gamma_{K3,2d}=\{\phi\in Aut(\Lambda_{K3})|<\phi(u),\phi(u)>=<u,u>$ and $\phi(u)=u\}.$
\end{center}

The moduli space of isomorphisms classes of polarized, algebraic K3 surfaces
of a fixed degree 2d, which we denote by $\mathcal{M}_{K3,pa}^{2d}$ is
isomorphic to a Zariski open set in the quasi-projective variety
$\Gamma_{K3,2d}$%
$\backslash$%
$\frak{h}_{2,19}.$

If we allow our surface to have singularities which are at most double
rational points, then the corresponding moduli space of polarized, algebraic
surfaces is equal to the entire locally symmetric space $\Gamma_{K3,2d}%
\backslash\frak{h}_{2,19}.$ In other words, marked and pseudo-polarized
surfaces corresponding to points in the complement of $\mathcal{M}%
_{K3,mpa}^{2d}$ in $\frak{h}_{2,19},$ are those surfaces for which the
projective image corresponding to any power of the polarization is singular
with singularities which are double rational points. Specifically, the
relation between $\mathcal{M}_{K3,mpa}^{2d}$and $\frak{h}_{2,19}$ is through
the period map, which we now will describe. The period map $\pi$ for marked K3
surfaces (X,$\alpha)$ is defined by integrating the holomorphic two form
$\omega$ along the basis $\alpha$ of $H_{2}(X,\mathbb{Z}),$ meaning

\begin{center}
$\pi(X,\alpha):=(...,\int_{\alpha_{i}}\omega,...)\in\mathbb{P}^{21}.$
\end{center}

The Riemann bilinear relations hold for $\pi(X,\alpha),$ meaning

\begin{center}
$\left\langle \pi(X,\alpha),\pi(X,\alpha)\right\rangle =0$ and \ $\left\langle
\pi(X,\alpha),\overline{\pi(X,\alpha)}\right\rangle >0.$
\end{center}

Choose a primitive vector $e\in H^{2}(X,\mathbb{Z)}$ such that $\left\langle
e,e\right\rangle =2d>0.$ As in \cite{PS}, one has the description of
$\frak{h}_{2,19}$ as one of the open sets of the quadric $Q$ defined in
$\mathbb{P}(\Lambda_{K3}\otimes\mathbb{C)}$ by the equations $\left\langle
u,u\right\rangle =0$ and $\left\langle u,e\right\rangle =0$ and the inequality
$\left\langle u,\overline{u}\right\rangle >0.$ Results from \cite{Ku} and
\cite{PS} combine to prove that the period map $\pi$ is surjection, i.e. each
point of the period domain $\frak{h}_{2,19}$ corresponds to a marked
pseudo-polarized algebraic K3 surface. By pseudo-polarized algebraic K3
surface we understand a pair (X,$e$) where e corresponds to either ample
divisor or pseudo ample divisor, which means that for any effective divisor
$D$ in X, we have $\left\langle D,e\right\rangle \geq0.$ Mayer proved the
linear system $|3e|$ defines a map:

\begin{center}
$\phi_{|3e|}:X\rightarrow X_{1}\subset\mathbb{P}^{m}$
\end{center}

such that: \textbf{i. }$X_{1}$ has singularities only double rational points.
\textbf{ii. }$\phi_{|3e|}$ is a holomorphic birational map. From the result of
Donaldson and the surjectivity of the period map, it follows that the moduli
space of pseudo-polarized algebraic K3 surfaces of degree 2d $\mathcal{M}%
_{K3,ppa}^{2d}$ is isomorphic to the locally symmetric space $\Gamma
_{K3,2d}\backslash\frak{h}_{2,19}.$ See \cite{D}. For the discussion of the
global Torelli Theorem for polarized, algebraic K3 surfaces see \cite{PS},
\cite{Ast} and \cite{BPV}.

\subsection{Description of the Discriminant Locus in the Moduli Space of
Pseudo-Polarized Algebraic K3 Surfaces}

\begin{notation}
From now on if a set $Y$ is contained in $X,$ then the complement of $Y$ in
$X$ will be denoted by $X\circleddash Y.$
\end{notation}

The complement of $\mathcal{M}_{K3,mpa}^{2d}$ in $\frak{h}_{2,19}$ can be
described as follow. Given a polarization class $e\in\Lambda_{K3},$ set
$T_{e}$ to be the orthogonal complement to $e$ in $\Lambda_{K3},$ i.e. $T_{e}
$ is the transcendental lattice. Then we have the realization of
$\frak{h}_{2,19}$ as one of the components of

\begin{center}
$\{u\in\mathbb{P}(T_{e}\otimes\mathbb{C})|\left\langle u,u\right\rangle =0$
and $\left\langle u,\overline{u}\right\rangle >0\}.$
\end{center}

Define the set $\Delta(e):=\{\delta\in\Lambda_{K3}|\left\langle e,\delta
\right\rangle =0$ and $\left\langle \delta,\delta\right\rangle =-2\}.$ For
each $\delta\in\Delta(e),$ define the hyperplane

\begin{center}
$H(\delta)=\{u\in\mathbb{P}(T_{e}\otimes\mathbb{C})|\left\langle
u,\delta\right\rangle =0\}.$
\end{center}

Let

\begin{center}
$\mathcal{H}_{K3,2d}=\underset{\delta\in\Delta(e)}{\cup}(H(\delta)\cap$
$\frak{h}_{2,19}).$
\end{center}

Set $\mathcal{D}_{K3}^{2d}:=\Gamma_{K3,2d}\backslash\mathcal{H}_{K3,2d}.$
Results from \cite{Ma}, \cite{PS}, \cite{To} and \cite{Ku} imply that
$\mathcal{D}_{K3}^{2d}$ is the complement of the moduli space of algebraic
polarized K3 surfaces $\mathcal{M}_{K3,pa}^{2d}$ in the locally symmetric
space $\Gamma_{K3,2d}\backslash\frak{h}_{K3,2d},$ i.e. $\mathcal{D}_{K3}%
^{2d}=(\Gamma_{K3,2d}\backslash\frak{h}_{K3,2d})\ \ominus\mathcal{M}%
_{K3,pa}^{2d}.$

\section{Enriques Surfaces and Their Moduli}

We will define an Enriques surface Y to be X/$\rho$, where X is a K3 surface
and $\rho$ is an involution acting on X without fixed points. On $\Lambda
_{K3}\approxeq\mathbf{H}^{3}\oplus(-E_{8})^{2}$ we will define the Enriques
involution $\rho(z_{1}\oplus z_{2}\oplus z_{3}\oplus x\oplus y)=(-z_{1}\oplus
z_{3}\oplus z_{2}\oplus y\oplus x).$ Let $\Lambda_{K3}^{+}$ and $\Lambda
_{K3}^{-}$ be the $\rho-$invariant and $\rho-$anti-invariant sublattices. The
unimodular lattice $\frac{1}{2}\Lambda_{K3}^{+}$ is isometric to the Enriques
lattice $\Lambda_{Enr}.$ We define the space $\Omega_{Enr}=\mathbf{P(}%
\Lambda_{K3}^{-}\otimes\mathbf{C)\cap}\Omega_{K3},$ where $\Omega_{K3}$ is the
period domain for marked K3 surfaces. It is easy to see that

\begin{center}
$\Omega_{Enr}:=SO_{0}(2,10)/SO(2)\times SO(10)=\frak{h}_{2,10}.$
\end{center}

\begin{definition}
\label{e1}We define $\Gamma_{Enr}$ $=rest_{\Lambda_{K3}^{-}}\{g\in
Aut(\Lambda_{K3})|g\circ\rho=\rho\circ g\}.$
\end{definition}

\begin{definition}
\label{e2}We defined $\widetilde{\Gamma}_{Enr}$ as a subgroup of finite index
in $\Gamma_{Enr}$ which preserves the so called $\mathbf{H-}$marking of the
Enriques surfaces, which means a pair (Y,j) and
\end{definition}

\begin{center}
j:$\mathbf{H\rightarrow}H^{2}(Y,\mathbb{Z})_{f}$.
\end{center}

\begin{remark}
\label{e3}It was proved in \cite{BPV} that $\widetilde{\Gamma}_{Enr}$ as a
subgroup of finite index in $\Gamma_{Enr}$.
\end{remark}

Let $\Delta_{+}:=\{\delta\in\Lambda_{K3}^{-}$
$\vert$%
$<\delta,\delta>=-2\}$ and $\Delta_{-}:=\{\delta\in\Lambda_{K3}|<\delta
,\delta>=-2$ and $\delta^{\rho}\neq\delta\}.$ It is shown on p. 283 of
\cite{BPV} that no point of the hyperplane $H_{l}=\{p\in\Omega_{Enr}%
|<p,\delta>=0$ and $\delta\in\Delta_{+}\}$ can be the period of a marked
Enriques surfaces. Namikawa showed in \cite{Na85} the following result. Let
$\delta\in\Delta_{-},$ then the points of $\ H_{\delta}=\{p\in\Omega
_{Enr}|<p,\delta-\delta^{\rho}>=0\}$ corresponds to Enriques surfaces with
double rational points$.$

\begin{definition}
\label{e5}We will define two divisors $\mathcal{D}_{+}:=\widetilde{\Gamma
}_{Enr}\backslash\cup_{l}H_{l}$ $\ ($for all $l\in\Delta_{-})$ and
$\mathcal{D}_{-}:=\widetilde{\Gamma}_{Enr}\backslash\cup_{\delta}H_{\delta}$
in $\mathcal{M}_{Enr,\mathbf{H}}^{2},$ (for all $\delta\in\Delta_{-}).$
\end{definition}

The following Theorems are due to Horikawa. See \cite{BPV}.

\begin{theorem}
The coarse moduli space of $\mathbf{H}-$marked Enriques surface $\mathcal{M}%
_{Enr,\mathbf{H}}^{2}$ is isomorphic to$\left(  \widetilde{\Gamma}%
_{Enr}\backslash\Omega_{Enr}\right)  \backslash\left(  \mathcal{D}_{+}%
\cup\mathcal{D}_{-}\right)  .$
\end{theorem}

\begin{theorem}
The isomorphism class of Enriques surface is uniquely determined by its period
in $\Gamma_{Enr}\backslash\Omega_{Enr}.$
\end{theorem}

\begin{remark}
\label{bor}The results of Borcherds in \cite{Bo96}, imply that both divisors
$\mathcal{D}_{+}$ and $\mathcal{D}_{-}$ are irreducible in $\mathcal{M}%
_{Enr,\mathbf{H}}^{2}$ and so their closures $\overline{\mathcal{D}_{+}}$ and
$\overline{\mathcal{D}_{-}}$ \ in the Baily-Borel compactification
$\overline{\widetilde{\Gamma}_{Enr}\backslash\Omega_{Enr}}=\overline
{\widetilde{\Gamma}_{Enr}\backslash\frak{h}_{2,10}}$ of $\ \widetilde{\Gamma
}_{Enr}\backslash\frak{h}_{2,10}.$
\end{remark}

\section{Discriminant as an Ample Divisor in the Moduli of Certain
Pseudo-Polarized Algebraic K3 Surfaces.}

We will assume from now on that the Baily-Borel compactification
$\overline{\Gamma_{K3,2d}\backslash\frak{h}_{K3,2d}}$ of $\Gamma
_{K3,2d}\backslash\frak{h}_{K3,2d}$ contains only one dimensional cusp. We
know that this is the case when d=1. See \cite{Hor}. Our proof that in this
case the closure $\overline{\mathcal{D}_{K3}^{2d}}$ of $\mathcal{D}_{K3}^{2d}$
in the Baily-Borel compactification $\overline{\Gamma_{K3,2d}\backslash
\frak{h}_{K3,2d}}$ of $\Gamma_{K3,2d}\backslash\frak{h}_{K3,2d}$ contains the
support of an ample divisor in $\overline{\Gamma_{K3,2d}\backslash
\frak{h}_{K3,2d}}$ will proceed in two steps. \textbf{Step 1.} We recall the
construction of an ample line bundle $\mathcal{L}_{K3,2d}^{-1}$ on
$\Gamma_{K3,2d}\backslash\frak{h}_{K3,2d}$ by means of the factor of
automorphy given by the functional determinant of the group action of
$\Gamma_{K3,2d}$ on $\frak{h}_{K3,2d}.$ Holomorphic sections of $\left(
\mathcal{L}_{K3,2d}^{-1}\right)  ^{\otimes n}$ are automorphic forms on
$\frak{h}_{K3,2d} $ with respect to the group $\Gamma_{K3,2d}.$ \textbf{Step
2.} We recall work from \cite{JT} which constructs an automorphic form on
$\frak{h}_{K3,2d} $ with respect to the group action $\Gamma_{K3,2d},$ and
which vanishes on a divisor $\mathcal{E}_{K3}^{2d}$ whose support is contained
in $\mathcal{D}_{K3}^{2d}.$

\subsection{Construction of a Line Bundle on $\frak{h}_{K3,2d}$}

Let us consider the principle bundle SO(2) bundle SO(2)$\rightarrow
$SO(2,19)/SO(19)$\rightarrow\frak{h}_{K3,2d}.$ Since SO(2)$=U(1)$ we can
associate a complex line bundle $\mathcal{L}_{\Gamma_{K3,2d}}\rightarrow
\Gamma_{K3,2d}\backslash\frak{h}_{K3,2d}.$ More explicitly, let $\gamma
\in\Gamma_{K3,2d}$ and set

\begin{center}
j($\gamma,\tau$):=$\det\left(  \frac{\partial\gamma(\tau)}{\partial\tau
}\right)  ^{-19}.$
\end{center}

Then $\mathcal{L}_{\Gamma_{K3,2d}}$ is constructed via the factor of
automorphy j($\gamma,\tau),$ i.e. $\mathcal{L}_{\Gamma_{K3,2d}}$ is obtained
as the quotient of $\frak{h}_{K3,2d}\times\mathbb{C}$ through the
identification $(\tau,w)\thicksim(\gamma\tau,j(\gamma,\tau)w).$ See
\cite{Tol}. From \cite{BB} and \cite{BB1}, we have that $\mathcal{L}%
_{\Gamma_{K3,2d}}^{-1}$ extends to an ample bundle on the Baily-Borel
compactification $\overline{\Gamma_{K3,2d}\backslash\frak{h}_{K3,2d}}$ of
$\Gamma_{K3,2d}\backslash\frak{h}_{K3,2d}$ , and automorphic forms on
$\frak{h}_{K3,2d}$ with respect to $\Gamma_{K3,2d}$ corresponds to sections of
$\left(  \overline{\mathcal{L}_{\Gamma_{K3,2d}}^{\ast}}\right)  ^{\otimes n},$
where $\overline{\mathcal{L}_{\Gamma_{K3,2d}}^{\ast}}$ is the extension of the
dual of the line bundle $\mathcal{L}_{\Gamma_{K3,2d}}$ to the Baily-Borel
compactification $\overline{\Gamma_{K3,2d}\backslash\frak{h}_{K3,2d}}$ of
$\mathcal{\Gamma}_{K3,2d}\backslash\frak{h}_{2,19}$ for some positive integer
$n.$ We have proved in \cite{JT} the following Lemma:

\begin{lemma}
\label{jt}Let $\pi:\mathcal{X\rightarrow M}_{K3,pa}^{2d}$ be the versal family
of polarized algebraic K3 surfaces of degree 2d, then $\pi_{\ast}%
(\mathcal{K}_{\mathcal{X}/\mathcal{M}_{K3,pa}^{2d}})\approxeq\mathcal{L}%
_{\Gamma_{K3,2d}}.$
\end{lemma}

The following Theorem from \cite{JT} is the main result of this section

\begin{theorem}
\label{jt1}Let $\overline{\Gamma_{K3,2d}\backslash\frak{h}_{K3,2d}}$ be the
Baily-Borel compactification of $\Gamma_{K3,2d}\backslash\frak{h}_{K3,2d}$
such that it contains only one zero dimensional cusp, then there is an ample
divisor $\overline{\mathcal{E}_{K3}^{2d}}$ on $\overline{\Gamma_{K3,2d}%
\backslash\frak{h}_{2,19}}$ with support contained in the closure of
$\mathcal{D}_{K3}^{2d}$ in $\overline{\Gamma_{K3,2d}\backslash\frak{h}%
_{K3,2d}}$ which we will denote by $\overline{\mathcal{D}_{K3}^{2d}}$ such
that $\mathcal{L}_{\Gamma_{K3,2d}}^{-1}\approxeq\mathcal{O}_{\Gamma
_{K3,2d}\backslash\frak{h}_{K3,2d}}(\mathcal{E}_{K3}^{2d}).$
\end{theorem}

\textbf{PROOF:} Our proof of Theorem \ref{jt} involves an explicit
constriction of a holomorphic function on $\frak{h}_{2,19}$ which is modular
with respect to $\Gamma_{K3,2d}$ and which is non-vanishing on $\mathcal{M}%
_{K3,pa}^{2d}$ and we know that $\mathcal{M}_{K3,pa}^{2d}=\left(
\Gamma_{K3,2d}\backslash\frak{h}_{K3,2d}\right)  \ominus\mathcal{D}_{K3}^{2d}.$

Since Baily-Borel bundle $\mathcal{L}_{\Gamma_{K3,2d}}^{-1}$ can be prolonged
to an ample bundle $\overline{\mathcal{L}_{\Gamma_{K3,2d}}^{-1}}$ over the
Baily-Borel compactification $\overline{\Gamma_{K3,2d}\backslash
\frak{h}_{2,19}}$, by the results in \cite{BB} and \cite{BB1} and since our
form is a section of a power of $\mathcal{L}_{\Gamma_{K3,2d}}^{-1},$ by Lemma
\ref{jt} hence Theorem \ref{jt1} follows.$\blacksquare$

Let (X,e) be a polarized K3 surface of degree 2d, and let $\mathcal{T}%
_{(X,e)}$ be the sheaf of holomorphic vector fields on (X,e). From
Kodaira-Spencer deformation theory, we can identify the tangent space
$T_{\mathcal{M}_{K3,mp}^{2d}}$ at the point (X,e) with $H^{1}(X,\mathcal{T}%
_{(X,e)}).$ The existence of the holomorphic two form $\omega$ on X implies
that we can identify $H^{1}(X,\mathcal{T}_{(X,e)})$ with $H^{1}%
(X,\mathcal{\Omega}^{1}),$ where $\Omega^{1}$ is the sheaf of holomorphic one
forms on X. One can deduce that the tangent space $T_{\mathcal{M}_{K3,mp}%
^{2d}}$ to the moduli space $\mathcal{M}_{K3,mp}^{2d}$ at the point
(X,$\alpha,e)$ can be identified with the space

\begin{center}
$H^{1}(X,\Omega^{1})_{0}=\{u\in H^{1}(X,\mathcal{\Omega}^{1})|<u,e>=0\}.$
\end{center}

We view any $\phi\in H^{1}(X,\mathcal{T}_{(X,e)})$ as a linear map from
$\Omega^{1,0}$ to $\Omega^{0,1}$ pointwise on X. Given $\phi_{1}$ and
$\phi_{2}$ in $H^{1}(X,\mathcal{T}_{(X,e)})$, the trace of the map $\phi
_{1}\circ\overline{\phi_{2}}:\Omega^{0,1}\rightarrow\Omega^{0,1}$ at a point a
point x$\in$X \ with respect to the unit volume CY metric g (meaning a
K\"{a}hler-Einstein metric compatible with the given polarization class e) is simply

\begin{center}
$Tr(\phi_{1}\circ\overline{\phi_{2}})=\underset{k,l,m,n}{\sum}\left(  \phi
_{1}\right)  _{\overline{l}}^{k}\overline{\left(  \phi_{2}\right)
_{\overline{n}}^{m}}g^{n\overline{l}}g_{k\overline{m}}.$
\end{center}

The existence of a Calabi-Yau metric on X compatible with the polarization e
(unique up to a scale) is guaranteed by Yau's Theorem \cite{Y}. We define
Weil-Petersson metric on $\mathcal{M}_{K3,mpa}^{2d}$ via the inner product

\begin{center}
$\left\langle \phi_{1},\phi_{2}\right\rangle :=\int_{X}Tr(\phi_{1}%
\circ\overline{\phi_{2}})vol_{g}.$
\end{center}

on the tangent space of $\mathcal{M}_{K3,mpa}^{2d}$ at (X,$\alpha,e)$. It is
shown in \cite{T1} that the Weil-Petersson metric is equal to the Bergman
metric on $\frak{h}_{2,19}.$ Therefore, the Weil-Petersson metric is a
K\"{a}hler metric with a K\"{a}hler form $\mu_{WP}.$

Since $\frak{h}_{2,19}$ is simply connected and over the moduli space of
marked polarized algebraic K3 surfaces $\mathcal{M}_{K3,mp}^{2d}%
\subset\frak{h}_{K3,2d}$ we have a universal family of marked and polarized K3
surfaces $\mathcal{X}^{2d}\rightarrow\mathcal{M}_{K3,mp}^{2d},$ there exists a
non-vanishing holomorphically varying family of holomorphic two forms over
$\frak{h}_{K3,2d}.$ For any such family consider the function on
$\frak{h}_{K3,2d}$ defined by

\begin{center}
$\left\|  \omega\right\|  _{L^{2}}^{2}=\left\langle \omega,\omega\right\rangle
=\int_{X}\omega\wedge\overline{\omega}.$
\end{center}

In \cite{T1} and \cite{Ti} it was proved that $\log\left\|  \omega\right\|
_{L^{2}}^{2}$ is a potential for the Weil-Petersson metric. The following
result from \cite{JT} proves the existence of a second potential for the
Weil-Petersson metric.

\begin{theorem}
\label{jt2}Let (X,e) be a polarized, algebraic K3 surface of degree 2d, and
let $\mu$ denote the unit volume Calabi-Yau form on X which is compatible with
the polarization class e. Let \{$\omega$\} be a non-vanishing, holomorphically
varying family of holomorphic two form on $\frak{h}_{2,19}.$ \textbf{A.} Let
$\det\Delta_{(X,e)}$ denote the zeta regularized product of the non-zero
eigenvalues of the Laplacian of the CY metric which acts on the space of
smooth functions on X. Then
\end{theorem}

\begin{center}
$dd^{c}\log\left(  \frac{\det\Delta_{(X,e})}{\left\|  \omega\right\|  _{L^{2}%
}^{2}}\right)  =0,$
\end{center}

\textit{or equivalently}

\begin{center}
$-dd^{c}\log\left(  \det\Delta_{(X,e})\right)  =-dd^{c}\log\left(  \left\|
\omega\right\|  _{L^{2}}^{2}\right)  =\mu_{WP}.$
\end{center}

\textit{in other words \ }$-\log\left(  \det\Delta_{(X,e})\right)  $
\textit{is a potential for the Weil-Petersson metric on }$\mathcal{M}%
_{K3,mpa}^{2d}.$

\textbf{B.}\textit{\ There is a holomorphic function (possibly multi-valued)
f}$_{K3,\omega,2d}$ on $\frak{h}_{2,19}\backslash\mathcal{M}_{K3,mpa}^{2d}$
\textit{such that}

\begin{center}%
$\vert$%
f$_{K3,\omega,2d}|^{2}=\left(  \frac{\det\Delta_{(X,e})}{\left\|
\omega\right\|  _{L^{2}}^{2}}\right)  ;$
\end{center}

\textit{hence} f$_{K3,\omega,2d}$ \textit{does not vanish on} $\mathcal{M}%
_{K3,mpa}^{2d}.$

The reader is referred to \cite{JT} for details of the proof of Theorem
\ref{jt2}.

\subsubsection{\textbf{Important Note}}

The function f$_{K3,\omega,2d}$ constructed in Theorem \ref{jt2} part
\textbf{B} is possibly multi-valued function with a divisor contained in the
complement of $\mathcal{M}_{K3,mpa}^{2d}$ in $\frak{h}_{2,19}.$ At this point,
we do not assert any behavior of f$_{K3,\omega,2d}$ with respect to the
discrete group $\Gamma_{K3,2d}.$ Theorem \ref{jt2} is valid for any degree 2d
of the polarization class $e$.

If we want to conclude automorphic behavior of the function f$_{K3,\omega,2d}$
we need to have that $\left\|  \omega\right\|  _{L^{2}}^{2}$ is a meromorphic
or holomorphic automorphic form since we know that $\det\Delta_{(X,e)}$ is a
function on $\mathcal{M}_{K3,pa}^{2d}=\Gamma_{K3,2d}\backslash\frak{h}_{2,19}.$

It is easy to see that we can always construct a meromorphic section
\{$\omega\}$ of the Baily-Borel line bundle $\mathcal{L}_{K3,2d}^{-1}.$ Then
f$_{K3,\omega,2d}$ will have additional zeroes and poles coming from the poles
and the zeroes of the meromorphic section \{$\omega\}$ of the Baily-Borel line
bundle $\mathcal{L}_{K3,2d}^{-1}.$

In order to conclude that some power of f$_{K3,\omega,2d}$ is an automorphic
form on $\Gamma_{K3,2d}\backslash\frak{h}_{2,19}$ we need to construct a
holomorphic family of non vanishing holomorphic two forms $\omega$ over
$\frak{h}_{2,19}$ such that $\left\|  \omega\right\|  _{L^{2}}^{2}$ is
automorphic form of weight -2. In order to construct such forms we will use
the special ,polarization 2d for which the Baily Borel compactification of the
moduli space of pseudo polarized K3 surface contain a unique cusp of dimension
0. The analogue of the non vanishing family of holomorphic forms in case of
elliptic curves is the family of the so called normalized one forms $\omega$
on the upper half plane $\frak{h}:=\{\tau\in\mathbb{C}|\operatorname{Im}%
\tau>0\}$ such that

\begin{center}
$\operatorname{Im}\tau=\left\|  \omega\right\|  _{L^{2}}^{2}=\frac{-\sqrt{-1}%
}{2}\int_{E}\omega\wedge\overline{\omega}.$
\end{center}

\begin{remark}
Theorem \ref{jt2} is a generalization of a known result which exists in the
setting of elliptic curves. A generalization of Theorem \ref{jt2} in the
setting of CY manifolds was established in \cite{JT1} and recently in
\cite{To99}.
\end{remark}

\subsection{Construction of a Special Family of Holomorphic Two Forms for Some Polarizations}

Our construction of the holomorphic family of normalized holomorphic two forms
in the case when the Baily-Borel compactification $\overline{\Gamma
_{K3,2d}\backslash\frak{H}_{2,19}}$ contains only one cusp of dimension 0.

One can construct a normalized family of forms $\{\omega_{2d}\}$ as follows:
\textbf{Step 1.} Let $\tau_{\infty}\in\overline{\Gamma_{K3,2d}\backslash
\frak{h}_{2,19}}$ is the unique zero dimensional cusp. Since $\overline
{\Gamma_{K3,2d}\backslash\frak{h}_{2,19}}$ is a projective variety, we can
find a disk D$\subset$ $\overline{\Gamma_{K3,2d}\backslash\frak{h}_{2,19}}$
such that $\tau_{\infty}\in D$ and D%
$\backslash$%
$\tau_{\infty}\subset$ $\Gamma_{K3,2d}\backslash\frak{h}_{2,19}.$ By
restricting the versal family $\pi_{2d}:\mathcal{X}^{2d}\rightarrow
\mathcal{M}_{K3,mp}^{2d}$ to the punctured disk D$^{\ast}$=D%
$\backslash$%
$\tau_{\infty}$ we obtain a family of algebraic K3 surfaces $\mathcal{X}%
_{D^{\ast}}\rightarrow D^{\ast}.$ The monodromy operator $T$ acting on
$H^{2}(X_{t},\mathbb{Z})$ is such that $(T^{m}-id)^{3}=0$ and $(T^{m}%
-id)^{2}\neq0.$ By taking a finite covering of $D^{\ast}$ we may assume that
$(T-id)^{3}=0$ and $(T-id)^{2}\neq0.$ \textbf{Step 2.} The family
$\mathcal{X}_{D^{\ast}}\rightarrow D^{\ast}$ constructed in \textbf{Step 1}
defines up to a sign a unique cycle $\gamma$ up to the action of the
automorphisms group of the primitive cycles such that $T\gamma=\gamma$ and
there exists cycles $\mu$ and $\eta$ such that $T\mu=\mu+\gamma$ and
$T\eta=\mu+\eta+\gamma.$ \textbf{Step 3.}Since $\frak{h}_{2,19}$ is a
contractible, there exists a globally defined, non-vanishing, holomorphically
varying family of holomorphic two forms, i.e.

\begin{center}
$\omega_{\tau}\in H^{0}(\frak{h}_{2,19},\pi_{\ast}\mathcal{K}_{\mathcal{X}%
_{K3}^{2d}/\mathcal{M}_{K3,mpp}^{2d}})$
\end{center}

where $\mathcal{K}_{\mathcal{X}_{K3}^{2d}/\mathcal{M}_{K3,mpp}^{2d}}$ is the
relative canonical sheaf. \textbf{Step 4.} The following Lemma is true:

\begin{lemma}
\label{ljt}The function $\phi(\tau):=\int_{\gamma}\omega_{\tau}$ is non
vanishing on $\frak{h}_{2,19}.$
\end{lemma}

\textbf{PROOF:} Suppose that at some point $\tau_{0}\in\frak{h}_{2,19}$
$\phi(\tau_{0})=0.$ From the epimorphism of the period map proved in \cite{To}
it follows that $\tau_{0}$ corresponds to the periods of some marked pseudo
polarized K3 surface (X$_{0},\alpha,e).$ It is easy to see that $\gamma\in
T_{e},$ i.e. $\left\langle \gamma,e\right\rangle =0.$ This implies that
$\gamma$ can be realized as an algebraic cycle in the K3 surface X$_{0}.$
Indeed we can find a line bundle $\mathcal{L}$ on X$_{0}$ such that the first
Chern class $c_{1}(\mathcal{L)}$ of $\mathcal{L}$ will be $\gamma.$ Then the
poles and zeroes of any meromorphic section of $\mathcal{L}$ will give a
realization of $\gamma$ as an algebraic cycle. A Theorem proved in \cite{PS}
states that after finite number of reflections generated by vectors $\delta
\in\Delta(e)\cap H^{1,1}(X_{0},\mathbb{R})$ we may assume that $\gamma$ can be
realized as an elliptic curve $E$ embedded in X$_{0}.$ On the other hand side
the condition $\left\langle \gamma,e\right\rangle =0$ implies that
$\left\langle E,E\right\rangle =\left\langle E,e\right\rangle =0.$ So from
here it follows that the Elliptic curve should be contracted by the map
$\phi_{|3e|}:X_{0}\rightarrow\mathbb{P}^{m}$ defined by the linear system
$|3e|.$ By a theorem of Grauert it is possible if and only if $\left\langle
E,E\right\rangle <0.$ See \cite{Gra}. So we obtain a contradiction. Lemma
\ref{ljt} is proved. $\blacksquare.$

\begin{definition}
\label{nor}We define the normalized family of holomorphic two forms as
\{$\omega_{n,2d}\}:=\left\{  \frac{\omega_{\tau}}{\phi(\tau)}\right\}  .$
\end{definition}

Following the identical steps in the case of elliptic curves, we construct a
family of forms $\omega$ such that

\begin{center}
$\operatorname{Im}\tau=\left\|  \omega\right\|  _{L^{2}}^{2}=\frac{-\sqrt{-1}%
}{2}\int_{E}\omega\wedge\overline{\omega}.$
\end{center}

\begin{definition}
\label{coc}For any g$\in\Gamma_{K3,2d}$ and K3 surface X represented by the
period point $\tau\in\frak{h}_{2,19}$, we set
\end{definition}

\begin{center}
$\psi(g,\tau)=\int_{\gamma}g^{\ast}\omega_{n,2d}$
\end{center}

\textit{where }$\gamma$\textit{\ denotes the invariant vanishing cycle.}

\begin{remark}
\label{coc1}$\psi(g,\tau)$ defines one cocycle with respect to the group
$\Gamma_{K3,2d}$ with coefficients in the $\Gamma_{K3,2d}$ module of
invertible analytic functions $\mathcal{O}_{\frak{h}_{2,19}}^{\ast}$ on
$\frak{h}_{2,19},$ i.e. \{$\psi(g,\tau)\}\in H^{1}(\Gamma_{K3,2d}%
,\mathcal{O}_{\frak{h}_{2,19}}^{\ast}).$
\end{remark}

\begin{lemma}
\label{coc2}When the degree of the polarization is such that there is a single
zero-dimensional cusp at the boundary of the Baily-Borel compactification
$\overline{\Gamma_{K3,2d}\backslash\frak{h}_{2,19}}$ of $\Gamma_{K3,2d}%
\backslash\frak{h}_{2,19},$ then the cocycle $\psi(g,\tau)$ defines a line
bundle on $\Gamma_{K3,2d}\backslash\frak{h}_{2,19}$ isomorphic to the
Baily-Borel line bundle $\mathcal{L}_{\Gamma_{K3,2d}}^{-1}.$
\end{lemma}

For detail proof of Lemma \ref{coc2} see Proposition \textbf{6.9.} of
\cite{JT}.

\begin{theorem}
\label{coc3}When the degree of the polarization is such that there is a single
zero-dimensional cusp at the boundary of the Baily-Borel compactification
$\overline{\Gamma_{K3,2d}\backslash\frak{h}_{2,19}}$ of $\Gamma_{K3,2d}%
\backslash\frak{h}_{2,19},$ then the function $\left\|  \omega_{n,2d}\right\|
_{L^{2}}^{2}$ on $\frak{h}_{2,19}$ is a modular form of weight $-2$ with
respect $\Gamma_{K3,2d}$ and it had no zeroes on $\frak{h}_{2,19}.$
\end{theorem}

\textbf{PROOF: }Theorem \ref{coc3} follows directly from Lemma \ref{coc2} and
the definition of the normalized holomorphic form. $\blacksquare.$

Now we can conclude that when the Baily-Borel compactification $\overline
{\Gamma_{K3,2d}\backslash\frak{H}_{2,19}}$ contains only one cusp of dimension
0 the function

\begin{center}
$|f_{K3,\omega_{n,2d}}|^{2}=\left(  \frac{\det\Delta_{(X,e})}{\left\|
\omega_{n,2de}\right\|  _{L^{2}}^{2}}\right)  $
\end{center}

is a modular function of weight $2$ with respect to $\Gamma_{K3,2d}$ and its
zero set is supported by $\mathcal{D}_{K3,2d}=\left(  \Gamma_{K3,2d}%
\backslash\frak{h}_{2,19}\right)  \setminus\mathcal{M}_{K3,pa}^{2d}. $ Notice
that $f_{K3,\omega_{n,2d}}$ is a holomorphic automorphic form defined up to a
character $\chi$ of the group $\Gamma_{K3,2d}/[\Gamma_{K3,2d},\Gamma
_{K3,2d}].$ A Theorem of Kazhdan states that this group is finite. See
\cite{Bour}. From here we conclude that ($f_{K3,\omega_{n,2d}})^{N},$ where
$N=\#\left(  \Gamma_{K3,2d}/[\Gamma_{K3,2d},\Gamma_{K3,2d}]\right)  $ will be
an automorphic form on $\Gamma_{K3,2d}\backslash\frak{h}_{2,19}$ whose zero
set is supported by $\mathcal{D}_{K3,2d}$ $,$ where $\mathcal{D}_{K3,2d}$ is
the complement of $\mathcal{M}_{K3,pa}^{2d}$ in $\left(  \Gamma_{K3,2d}%
\backslash\frak{h}_{2,19}\right)  ,$ i.e. $\mathcal{D}_{K3,2d}=\left(
\Gamma_{K3,2d}\backslash\frak{h}_{2,19}\right)  \circleddash\mathcal{M}%
_{K3,pa}^{2d}.$ Thus Theorem \ref{jt1} is proved. $\blacksquare.$

We note that the case d=2 is the main point of consideration in \cite{BKPS}.
In \cite{BKPS} Borcherds constructed an automorphic form for degree two
polarization K3 surfaces, whose zero set is supported by the discriminant locus.

\section{Uniform Bounds.}

\subsection{Uniform Bounds for K3 Surfaces}

Given a family of polarized, algebraic K3 or Enriques surfaces fibred over a
curve C, we can give a short argument determining the number of singular fibres.

\begin{theorem}
\label{b} Let C be an algebraic curve of genus g and $\pi:Y\rightarrow C$ be a
three dimensional projective non-singular variety such that for every $t\in
C,$ $\pi^{-1}(t)=X_{t}$ is a non-singular K3 surface. Suppose that on $Y$ we
have a polarization class $H$ such that $H\left|  _{X_{t}}\right.  =e$ and
$<$%
e,e%
$>$%
=2d, where d is any positive integer. Let m$_{\infty}$ be the number of points
on C for which the local monodromy operator is of infinite order. Then the
number of singular fibres of $\pi$ is less or equal to $2g-2+m_{\infty}.$
\end{theorem}

\textbf{PROOF:} Consider the relative dualizing line bundle $\mathcal{K}%
_{\mathcal{X}/\mathcal{M}_{K3,ppa}^{2d}}$ of the universal family of marked
polarized algebraic

\begin{center}
$\pi:\mathcal{X}^{2d}\mathcal{\rightarrow M}_{K3,mpa}^{2d}=\frak{h}_{2,19}$.
\end{center}

On $\pi_{\ast}\mathcal{K}_{\mathcal{X}/\mathcal{M}_{K3,mpa}^{2d}}$ we have a
natural metric. Indeed the fibres of the line bundle $\pi_{\ast}%
\mathcal{K}_{\mathcal{X}/\mathcal{M}_{K3,mpa}^{2d}}$ over $\tau\in
\mathcal{M}_{K3,mpa}^{2d}$ is $H^{0}(X_{\tau},\Omega_{\Xi_{\tau}}^{2}).$ So
the natural metric will be

\begin{center}
$\left\|  \cdot\right\|  ^{2}=\left\langle \omega_{\tau},\omega_{\tau
}\right\rangle =\int_{X_{\tau}}\omega_{\tau}\wedge\overline{\omega_{\tau}}.$
\end{center}

In \cite{To} it is shown that one has the formula $c_{1}(\left\|
\cdot\right\|  ^{2})=-\mu_{B}$ where $\mu_{B}$ is the form associated with the
complete Bergman metric on $\frak{h}_{2,19}.$ Hence, $-c_{1}(\left\|
\cdot\right\|  )$ defines a compete metric on $C$
$\backslash$%
$E_{\infty}$ where $E_{\infty}$ is the set of points of $C$ around which the
local monodromy is infinite. Note that $-c_{1}(\left\|  \cdot\right\|  )$ when
restricted to $C$
$\backslash$%
$E_{\infty}$ is integrable. This is so since the periods

\begin{center}
$\left(  ...,\int_{\gamma_{i}}\omega_{t},..\right)  $
\end{center}

are solutions of ordinary differential equations with regular solutions, hence
have logarithmic growth near $E_{\infty}$ and

\begin{center}
$\left\langle \omega_{\tau},\omega_{\tau}\right\rangle =\left(  ...,\int
_{\gamma_{i}}\omega_{t},..\right)  \left(  \left\langle \gamma_{i},\gamma
_{j}\right\rangle \right)  \left(  ...,\overline{\int_{\gamma_{i}}\omega_{t}%
},..\right)  ^{t}.$
\end{center}

By the Gauss-Bonnet theorem, we have

\begin{center}
$\int_{C\backslash E_{\infty}}-c_{1}(\left\|  \cdot\right\|  )=-\chi(C$
$\backslash$%
$E_{\infty})=2g-2+m_{\infty}.$
\end{center}

From Borcherds's result proved in \cite{BKPS}, that one can always construct a
holomorphic automorphic form for any polarization class e on $\Gamma
_{2d}\backslash\frak{h}_{2,19}$, whose zero set $\mathcal{H}_{K3,2d}$ contains
the support of the discriminant locus $\mathcal{D}_{K3}^{2d},$ we observe that
the number of points on $C$ corresponding to singular fibres is less or equal
to $-\chi(C\backslash E)$ since

\begin{center}
$\left\langle \mathcal{H}_{K3,2d},C\right\rangle =\int_{C\backslash E_{\infty
}}-c_{1}(\left\|  \cdot\right\|  )=-\chi(C\backslash E).$
\end{center}

This proves Theorem \ref{b}. $\blacksquare.$

\subsection{Uniform Bounds for Enriques Surfaces}

\begin{theorem}
\label{b1}Let $\mathcal{Y\rightarrow}C$ be a family of $\mathbf{H}$ marked
Enriques surface over the algebraic curve C. (The $\mathbf{H}$ marking is
defined in Definition \ref{e2}.) Let m$_{\infty}$ be the number of points on C
for which the local monodromy operator is of infinite order. Then the number
of singular fibres of $\pi$ is less or equal to $2(2g-2+m_{\infty}). $
\end{theorem}

\textbf{PROOF:} The proof of Theorem \ref{b1} is the same as the proof of
Theorem \ref{b} by taking into account Remark \ref{bor} which states that the
discriminant locus in the moduli space $\widetilde{\Gamma}_{Enr}%
\backslash\Omega_{Enr}=\mathcal{M}_{Enr,\mathbf{H}}$ consists of
$\mathcal{D}_{+}$ and $\mathcal{D}_{-}$ and $\mathcal{L}_{\widetilde{\Gamma
}_{Enr}}^{-1}=\chi^{+}\otimes\mathcal{O(D}_{+})$ and $\mathcal{L}%
_{\widetilde{\Gamma}_{Enr}}^{-1}=\chi^{-}\otimes\mathcal{O(D}_{-}).$ Here
$\chi^{+}$ and $\chi^{-}$ are characters of the finite group $\widetilde
{\Gamma}_{Enr}/[\widetilde{\Gamma}_{Enr},\widetilde{\Gamma}_{Enr}].$ Thus
Theorem \ref{b1} is proved. $\blacksquare.$

\section{Finiteness Theorems.}

\subsection{K3 Surfaces}

Given the structural results from \textbf{Section 2} and the construction of
an automorphic form from \textbf{Section 4,} we now can prove finiteness for
certain families of K3 surfaces. As before, the K3 surface is assumed to have
a polarization of degree 2d such that the Baily-Borel compactification
$\overline{\Gamma_{K3,2d}\backslash\frak{H}_{2,19}}$ contains only one cusp of
dimension 0. We know from the result of Horikawa that for degree 2 K3 surfaces
the Baily-Borel compactification $\overline{\Gamma_{K3,2}\backslash
\frak{H}_{2,19}}$ contains only one cusp of dimension 0.

\begin{definition}
\label{D1}Let C be a fixed, non-singular algebraic curve, and let
E$\neq\emptyset$ and E be a fixed effective divisor on C such that all the
points in E have multiplicity 1. Define Sh(C,E,\{K3,2d\}) to be the set of all
isomorphism classes of three dimensional projective, polarized varieties which
admit fibration over C such that any fibre over C
$\backslash$%
E is a non-singular K3 surface for which the induced polarization is such that
the Baily Borel compactification of the coarse moduli space contains only one
cusp of dimension zero.
\end{definition}

\begin{theorem}
\label{Th1}Let C be an algebraic curve and let $\pi:Y\rightarrow C$ be a three
dimensional projective non-singular variety such that for every $t\in C,$
$\pi^{-1}(t)=X_{t}$ is a non-singular K3 surface. Suppose that on $Y$ we have
a polarization class $H$ such that $H\left|  _{X_{t}}\right.  =e$ and
$<$%
e,e%
$>$%
=2d is such that the Baily-Borel compactification $\overline{\Gamma
_{K3,2d}\backslash\frak{H}_{2,19}}$ contains only one cusp of dimension 0.
Then the family $\pi:Y\rightarrow C$ is isotrivial.
\end{theorem}

\textbf{PROOF:} From the versal properties of the moduli space $\mathcal{M}%
_{K3,ppa\text{ }}^{2d}$ (See \cite{PS}) it follows that we have a map

\begin{center}
$p:C\rightarrow\mathcal{M}_{K3,ppa\text{ }}^{2d}=\Gamma_{2d}\backslash
\frak{h}_{2,19}\subset\overline{\mathcal{M}_{K3,ppa\text{ }}^{2d}}$
\end{center}

such that $p(C)\cap\mathcal{D}_{K3}^{2d}=\emptyset.$ So we deduce that $p(C)$
is contained in the complement of $\overline{\mathcal{D}_{K3}^{2d}}$ in
$\overline{\mathcal{M}_{K3,ppa\text{ }}^{2d}},$ where $\overline
{\mathcal{M}_{K3,ppa\text{ }}^{2d}}$ is the Baily-Borel compactification of
$\mathcal{M}_{K3,ppa\text{ }}^{2d}=\Gamma_{K3,2d}\backslash\frak{h}_{2,19}$
and $\overline{\mathcal{D}_{K3}^{2d}}$ is the closure of $\mathcal{D}%
_{K3}^{2d}$ in $\overline{\mathcal{M}_{K3,ppa}^{2d}}.$ Since $\overline
{\mathcal{D}_{K3}^{2d}}$ is an ample divisor, we deduced that the complement
of $\overline{\mathcal{D}_{K3}^{2d}}$ in $\overline{\mathcal{M}_{K3,ppa}^{2d}%
}$ is an affine variety. We deduce that $p(C)$ must be a point, since $p(C)$
is a projective variety in the affine variety $\overline{\mathcal{M}%
_{K3,ppa\text{ }}^{2d}}$ $\circleddash$ $\overline{\mathcal{D}_{K3}^{2d}}$ .
This proved Theorem \ref{Th1}. $\blacksquare.$

\begin{theorem}
\label{Th2} Suppose that we fix an algebraic curve C, a divisor E$\neq
\emptyset$ on C as in Definition \ref{D1} and degree of polarization
$<$%
e,e%
$>$%
=2d such that the Baily-Borel compactification $\overline{\Gamma
_{K3,2d}\backslash\frak{H}_{2,19}}$ contains only one cusp of dimension 0,
then the set Sh(C,E,\{K3,2d\}) is finite.
\end{theorem}

\textbf{PROOF:} The proof will be done in two steps. \textbf{Step1. }We will
prove that there exists a holomorphic map p:$C\rightarrow\overline
{\Gamma_{K3,2d}\backslash\frak{h}_{2,19}}=\overline{\mathcal{M}_{K3,ppa\text{
}}^{2d}}.$ \textbf{Step2.} We will apply Theorem \ref{t1} to the pair $C$ and
$\overline{\mathcal{M}_{K3,ppa\text{ }}^{2d}}$ and the ample divisor
$\overline{\mathcal{D}_{K3}^{2d}}$ $\subset$ $\overline{\mathcal{M}%
_{K3,ppa\text{ }}^{2d}}$ taking into account that by Theorem \ref{b} p($C)$
have a bounded volume to conclude Theorem \ref{Th2}.

\begin{proposition}
\label{Th21}There exists a holomorphic map p:$C\rightarrow\overline
{\Gamma_{K3,2d}\backslash\frak{h}_{2,19}}.$ \textit{such that }$\overline
{p}^{\ast}\overline{\mathcal{D}_{K3}^{2d}}=E.$
\end{proposition}

\textbf{PROOF: }From the versal properties of the moduli space $\Gamma
_{K3,2d}\backslash\frak{h}_{2,19}$ of pseudo-polarized algebraic K3 surfaces
it follows that we have a map

\begin{center}
$p:C$ $\backslash$ $E\rightarrow\Gamma_{K3,2d}\backslash\frak{h}_{2,19}%
\subset\overline{\Gamma_{K3,2d}\backslash\frak{h}_{2,19}}=\overline
{\mathcal{M}_{K3,ppa\text{ }}^{2d}}$
\end{center}

such that $p(C$
$\backslash$%
$E)\cap\mathcal{D}_{K3}^{2d}=\emptyset.$ (See \cite{PS}.) Let us denote by
$E_{f}$ those points around which the monodromy have a finite order and by
$E_{\infty}$ those points around which the monodromy is of infinite order.
From Theorem 9.5. of Ph. Griffith proved in \cite{Gr} we can conclude that the
map $p$ can be prolonged through C$\backslash E_{\infty},$ i.e. we have $p:C$
$\backslash$ $E_{\infty}\rightarrow\Gamma_{K3,2d}\backslash\frak{h}_{2,19}$
and p($E_{f})\subset\mathcal{D}_{K3}^{2d}.$ Borel proved in \cite{B}:

\begin{theorem}
\label{B}Let $\mathcal{S}$ be a bounded, symmetric domain and $\Gamma\subset
Aut(\mathcal{S})$ an arithmetically defined discrete group of automorphisms.
Let $\mathcal{U}:=\Gamma\backslash\mathcal{S}.$ Let p:D$^{\ast}\rightarrow
\mathcal{U}$ be a holomorphic map from the punctured Disk D$^{\ast}$ to the
locally symmetric space $\mathcal{U}$, then p can be extended to a holomorphic
map $\overline{p}:D\rightarrow\overline{\mathcal{U}},$ where $\overline
{\mathcal{U}}$ is the Baily-Borel compactification of $\mathcal{U}.$
\end{theorem}

Theorem \ref{B} implies that there is a holomorphic map $\overline
{p}:C\rightarrow\overline{\Gamma_{K3,2d}\backslash\frak{h}_{2,19}}$ such that
$Supp\left(  (\overline{p})^{\ast}(\overline{\mathcal{D}_{K3}^{2d}})\right)
=E.$ Proposition \ref{Th21} is proved. $\blacksquare.$

For the curve C the map $\overline{p}:C\rightarrow\overline{\Gamma
_{K3,2d}\backslash\frak{h}_{2,19}}$ will be either finite map onto
$\overline{p}(C)$ or it will be a map to a point. The last possibility is
impossible since this will mean that the family $\pi:Y\rightarrow C$ is
isotrivial. We assumed that this is not the case. From here and Theorem
\ref{b} we can see that conditions \textbf{a, b }and \textbf{c }of Theorem
\ref{t1} are satisfied so we conclude that Sh(S,E,\{K3,2d\}) is finite.
Theorem \ref{Th2} is proved. $\blacksquare.$

\subsection{Finiteness Theorems for S-K3 Surfaces}

Following \cite{BKPS} we define an S-K3 surface X for some Lorenzian lattice
S$\subset\Lambda_{K3}$ of rank $\geq1$ is a K3 surface with a fixed primitive
embedding of S into the Picard group such that the image of S contains a
pseudo-ample class. (A pseudo-ample class is a class D such that $\left\langle
D,D\right\rangle >0$ and $\left\langle D,C\right\rangle \geq0$ for all curves
C on the K3 surface X.

It follows from the surjectivity of the period map that the moduli space of
S-K3 surfaces is isomorphic to $\Gamma_{S}\backslash\frak{h}_{2,20-rkS},$
where $\Gamma_{S}$ is an arithmetic group acting on $\frak{h}_{2,20-rkS}.$ See
\cite{T1}. We can define in $\Gamma_{S}\backslash\frak{h}_{2,20-rkS}$ the
discriminant locus as in the case of polarized algebraic K3 surfaces. We
define $T_{S}$ as follows: $T_{S}:=\{u\in\Lambda_{K3}|\left\langle
u,S\right\rangle =0\}.$ Define the set $\Delta(S):\Delta(S):=\{\delta\in
T_{S}|\left\langle \delta,\delta\right\rangle =-2\}.$ Remember that
$\frak{h}_{2,20-rkS}$ is one of the component of the set $\{u\in
\mathbb{P}(T_{S}\otimes\mathbb{C})|\left\langle u,u\right\rangle
=0\&\left\langle u,u\right\rangle >0\}.$ For each $\delta\in\Delta_{S},$
define the hyperplane $H(\delta)=\{u\in\mathbb{P}(T_{e}\otimes\mathbb{C}%
)|\left\langle u,\delta\right\rangle =0\}.$ Let $\mathcal{H}_{S}%
=\underset{\delta\in\Delta(S)}{\cup}(H(\delta)\cap$ $\frak{h}_{2,20-rkS}).$
Set $\mathcal{D}_{S}:=\Gamma_{S}\backslash\mathcal{H}_{S}.$ \ 

\begin{definition}
\label{DD1}Let us define Sh(C,E,(K3,S)) as the set of \ families of algebraic
S-K3 surfaces up to isomorphisms over the algebraic curve C such that for each
t$\in$C%
$\backslash$%
E, $\pi^{-1}(t)$ is a non-singular S-K3 surface and for each t$\in$E,
$\pi^{-1}(t)$ is a singular surface.
\end{definition}

Applying the same arguments as in the previous Section we obtain the following Theorem:

\begin{theorem}
Suppose that we fix an algebraic curve C, a divisor E$\neq\emptyset$ on C as
in Definition \ref{DD1} and suppose that on $\frak{h}_{2,20-rkS}$ there exists
a holomorphic automorphic $\Phi_{S}$ form such that the supporter of the zero
set of $\Phi_{S}$ is exactly $\mathcal{D}_{S},$ then the set Sh(C,E,\{K3,S\})
is finite.
\end{theorem}

Borcherds and Nikulin found some lattices S as NS groups of K3 surfaces for
which one can construct a holomorphic automorphic $\Phi_{S}$ form such that
the supporter of the zero set of $\Phi_{S}$ is exactly $\mathcal{H}_{S}.$ See
\cite{BKPS}.

\subsection{Finiteness Theorems for Enriques Surfaces}

\begin{definition}
\label{D2}Let C be a fixed, non-singular algebraic curve, and let E be a fixed
effective divisor on C such that all the points in E have multiplicity 1.
Define Sh(C,E,Enr) to be the set of all isomorphism classes of three
dimensional projective, polarized varieties which admit fibration over C such
that the fibres over C
$\backslash$%
E are non-singular Enriques surfaces for which the induced polarization is of
degree 2d \ for d%
$>$%
1 or they are $\mathbf{H}$ polarized.
\end{definition}

\begin{theorem}
\label{Enr}Let C be an algebraic curve and let $\pi:Y\rightarrow C$ be a three
dimensional projective non-singular variety such that for every $t\in C,$
$\pi^{-1}(t)=X_{t}$ is a non-singular Enriques surface. Suppose that on $Y$ we
have a polarization class $H$ such that $H\left|  _{X_{t}}\right.  =e$ and
$<$%
e,e%
$>$%
=2d for d%
$>$%
1 or all the fibres are $\mathbf{H}$ marked. Then the family $\pi:Y\rightarrow
C$ is isotrivial.
\end{theorem}

\textbf{PROOF: }We will define the groups $\Gamma_{Enr,2d}$ for $d>1$ as follows:

\begin{definition}
\label{e4}Let $e\in\Lambda_{K3}$ be a primitive element such that $<e,e>=4d $
for $d>1,$ then
\end{definition}

\begin{center}
$\Gamma_{Enr,2d}:=rest_{\Lambda_{K3}^{-}}\{g\in Aut(\Lambda_{K3})|g\circ
\rho=\rho\circ g$ and $g(e)=e\}.$
\end{center}

Clearly $\Gamma_{Enr,2d}$ is a subgroup of finite index in $\Gamma_{Enr}$ as
defined in Definition \ref{e1}.

First we will define the moduli space $\mathcal{M}_{Enr,2d}$ of polarized
Enriques surfaces with degree of polarization 2d%
$>$%
2. We already defined

\begin{center}
$\mathcal{M}_{Enr,\mathbf{H}}=\left(  \widetilde{\Gamma}_{Enr}\backslash
\frak{h}_{2,10}\right)  \circleddash\left(  \mathcal{D}_{+})\cup
\mathcal{D}_{-}\right)  .$
\end{center}

It was proved in \cite{BPV} that $\mathcal{M}_{Enr,\mathbf{H}}$ is the coarse
moduli space of $\mathbf{H}$ marked Enriques surfaces.

\begin{definition}
We will define $\mathcal{M}_{Enr,2d}$ as follows:
\end{definition}

\begin{center}
$\mathcal{M}_{Enr,2d}:=\left(  \Gamma_{Enr,2d}\backslash\frak{h}%
_{2,10}\right)  \circleddash(\pi_{2d}^{-1}(\mathcal{D}_{E,+})\cup\pi_{2d}%
^{-1}(\mathcal{D}_{E,-})),$
\end{center}

\textit{where }$\pi_{2d}:\Gamma_{Enr,2d}\backslash\frak{h}_{2,10}%
\rightarrow\Gamma_{Enr}\backslash\frak{h}_{2,10}$ \textit{is the natural map
and }$\pi_{\mathbf{H}}:\widetilde{\Gamma}_{Enr}\backslash\frak{h}%
_{2,10}\rightarrow\Gamma_{Enr}\backslash\frak{h}_{2,10}.$ $\mathcal{D}_{\pm}$
\textit{are defined in Definition} \ref{e5} and $\mathcal{D}_{E,\pm}=\left(
\pi_{\mathbf{H}}\right)  _{\ast}(\mathcal{D}_{\pm}).$

Global Torelli Theorem for Enriques surfaces implies that $\mathcal{M}%
_{Enr,2d}$ is coarse moduli space of Enriques surfaces with 2d polarization.
It follows from the versal properties of the coarse moduli space
$\mathcal{M}_{Enr,2d}$ of polarized Enriques surfaces for $d>1$ that we have a
holomorphic map $p:C\rightarrow\mathcal{M}_{Enr,2d}.$ Since $\pi_{\mathbf{H}%
}:\widetilde{\Gamma}_{Enr}\backslash\frak{h}_{2,10}\rightarrow\Gamma
_{Enr}\backslash\frak{h}_{2,10}$ is a finite map, then the image of an affine
open set or quasi affine are also affine or quasi-affine. We know that
according to \cite{Bo96}

\begin{center}
$\left(  \widetilde{\Gamma}_{Enr}\backslash\frak{h}_{2,10}\right)
\ominus\left(  \mathcal{D}_{+})\cup\mathcal{D}_{-}\right)  $
\end{center}

is quasi-affine. From here we obtain that $\mathcal{M}_{Enr,2d}$ is
quasi-affine. So the map $p:C\rightarrow\mathcal{M}_{Enr,2d}$ is a map to a
point. The same arguments are applied when we consider $\mathbf{H}$ marked
Enriques surfaces. Thus we proved Theorem \ref{Enr}. $\blacksquare.$

\begin{theorem}
\label{B1}Suppose that we fix and algebraic curve C, a divisor E$\neq
\emptyset$ as in Definition \ref{D2}, and degree of polarization
$<$%
e,e%
$>$%
=2d$.$ Then the set Sh(C,E,Enr) is finite.
\end{theorem}

\textbf{PROOF:} The proof of Theorem \ref{B1} is reduced to the case of
$\mathbf{H}$ marked Enriques surfaces and then the proof goes exactly in the
same way as in the case of pseudo polarized algebraic K3 surfaces whose moduli
space has a unique cusp of dimension zero in the Baily Borel compactification.

From the versal properties of the coarse moduli space of polarized Enriques
surfaces $\mathcal{M}_{Enr,2d}$ and by using the results of Griffiths and
Borel and the properties of the Baily-Borel compactification, we get a
holomorphic map:

\begin{center}
$p:C\rightarrow\overline{\Gamma_{Enr,2d}\backslash\frak{h}_{2,19}}$
\end{center}

where $\overline{\Gamma_{Enr,2d}\backslash\frak{h}_{2,19}}$ is the Baily-Borel
compactification \ of $\Gamma_{Enr,2d}\backslash\frak{h}_{2,19}.$ Since
$\Gamma_{Enr,2d}$ is a subgroup of finite index in $\Gamma_{Enr}$ we obtain a
holomorphic map $p_{1}:C\rightarrow\overline{\Gamma_{Enr}\backslash
\frak{h}_{2,19}}$. We have a finite\ fixed map according to \cite{BB} and
\cite{BB1}: $\pi_{d}:\overline{\Gamma_{Enr,2d}\backslash\frak{h}_{2,19}%
}\rightarrow\overline{\Gamma_{Enr}\backslash\frak{h}_{2,19}}.$ We also know
that $\widetilde{\Gamma}_{Enr}$ is a subgroup of finite index in $\Gamma
_{Enr}.$ So we have a finite map: $\pi_{\mathbf{H}}:\overline{\widetilde
{\Gamma}_{Enr}\backslash\frak{h}_{2,19}}\rightarrow\overline{\Gamma
_{Enr}\backslash\frak{h}_{2,19}}.$

After taking a finite cover

\begin{center}
$\psi:\widetilde{C}=C\times_{p_{1}(C)}\pi_{\mathbf{H}}^{-1}(p_{1}%
(C))\rightarrow C$
\end{center}

of a degree less or equal to $\deg\pi_{\mathbf{H}},$ we will get a finite
holomorphic map

\begin{center}
$\widetilde{p}:\widetilde{C}\rightarrow\widetilde{p}(\widetilde{C}%
)\subset\overline{\widetilde{\Gamma}_{Enr}\backslash\frak{h}_{2,19}}.$
\end{center}

Over $\widetilde{C}$ we have a family of Enriques surfaces with $\mathbf{H}$
polaization and discriminant locus $\widetilde{E}=\psi^{-1}(E).$ We can apply
now Theorem \ref{b1} to get a bound on the number of points in $\widetilde{E}$
and the volume of the image of $\widetilde{p}(C)$ in $\overline{\widetilde
{\Gamma}_{Enr}\backslash\frak{h}_{2,19}}.$ Since the degree of the map $\psi$
is fixed we get a bound on the volume of the image $p_{1}(C)$ in
$\overline{\Gamma_{Enr}\backslash\frak{h}_{2,19}}$ and respectively of the
volume of the image $p(C)$ in $\overline{\Gamma_{Enr,2d}\backslash
\frak{h}_{2,19}}.$ Notice also that from Remark \ref{bor} we conclude that
$(\pi_{d})^{\ast}(\pi_{\mathbf{H}})_{\ast}(\mathcal{D}_{+}+\mathcal{D}_{-})$
is an ample divisor. Since $p(E)\subset(\pi_{d})^{\ast}(\pi_{\mathbf{H}%
})_{\ast}(\mathcal{D}_{+}+\mathcal{D}_{-}),$ Theorem \ref{t1} implies that the
set Sh($C,E,Enr)$ is finite by repeating the arguments of Theorem \ref{Th2}.
From here Theorem \ref{B1} follows directly. $\blacksquare.$

\section{Isotriviality.}

\begin{definition}
We will say that a family of algebraic varieties $X\rightarrow Y$ is an
isotrivial family if there exists a finite map $\phi:Y_{1}\rightarrow Y$ such
that the family $X\times_{Y}Y_{1}\rightarrow Y_{1}$ is a trivial one, i.e. the
family $X\times_{Y}Y_{1}\rightarrow Y_{1}$ is isomorphic to $Y_{1}\times Z.$
\end{definition}

\begin{theorem}
Suppose that $\mathcal{M}$ is the coarse moduli space of polarized algebraic
varieties $Z$. Suppose that $\mathcal{M}$ is a quasi projective variety. Let
$\overline{\mathcal{M}}$ be some projective compactification of $\mathcal{M}$
such that $\overline{\mathcal{M}}$ $\ominus\mathcal{M=D}$ is a divisor with
normal crossings. If $\mathcal{D}$ supports an ample divisor, then any family
$\mathcal{Z}\rightarrow C$ of algebraic polarized varieties $Z$ over a
projective variety $C$ without singular fibres is isotrivial.
\end{theorem}

\textbf{PROOF:} The proof is obvious since the condition that $\mathcal{D}$
supports an ample divisor implies that $\mathcal{M}$ is an affine or
quasi-affine. On the other hand side from the versal properties of the course
moduli space we deduce that there is a map $p:C\rightarrow\mathcal{M}.$ Since
$C$ is a projective variety, then $p(C)$ must be a point. $\blacksquare.$ See
also \cite{BKPS}.

\end{document}